\documentclass[12pt]{article}

\usepackage{latexsym}
\usepackage{amsopn,amscd}
\usepackage{amsfonts}
\usepackage{amssymb}
\usepackage{amsthm}
\usepackage{amsmath}
\usepackage{a4wide}

\newtheorem{Theorem}{Theorem}[section]
\newtheorem{Proposition}[Theorem]{Proposition}
\newtheorem{Lemma}[Theorem]{Lemma}

\newtheorem{Remark}[Theorem]{Remark}

\def\bemphas{}

\def\scs{\sc}
\def\dd{d}

\def\ppartial{\mathcal D}
\def\pppartial{\mathcal D}
\def\ppppartial{[\, \cdot \, , \, \cdot \, ]}

\def\nnabla{\widetilde \nabla}
\def\HH{\widetilde H}

\long
\def\MSC#1\EndMSC{\def\arg{#1}\ifx\arg\empty\relax\else
      {\par\narrower\noindent
      2000 Mathematics Subject Classification: #1\par}\fi}

\long
\def\KEY#1\EndKEY{\def\arg{#1}\ifx\arg\empty\relax\else
    {\par\narrower\noindent
      Keywords and Phrases: #1\par}\fi}
\def\Nsddata#1#2#3#4#5{
\begin{equation*}
   (#4
     \begin{CD}
      \null @>#2>> \null\\[-3.2ex]
      \null @<<#3< \null
     \end{CD}
    #1, #5)
  \end{equation*}
}

\def\fra{\mathfrak}
\def\Sigm{\mathcal S}

\numberwithin{equation}{section}

\begin{document}

\title{\bf The sh-Lie algebra perturbation Lemma}

\author{
J.~Huebschmann
\\[0.3cm]
 USTL, UFR de Math\'ematiques\\
CNRS-UMR 8524
\\
59655 Villeneuve d'Ascq Cedex, France\\
Johannes.Huebschmann@math.univ-lille1.fr
 }
\maketitle

\begin{abstract}
{Let $R$ be a commutative ring with $1$ which contains
the rationals as a subring and let   $(M
     \begin{CD}
      \null @>{\nabla}>> \null\\[-3.2ex]
      \null @<<{\pi}< \null
     \end{CD}
    \fra g, h)$
be a {\em contraction\/} of chain complexes (over $R$). We denote the symmetric coalgebra functor by
$\Sigm^{\mathrm c}$, the loop Lie algebra functor by $\mathcal L$,
the classifying coalgebra functor by $\mathcal C$, and the suspension operator
by $s$. We shall
establish the following.

\noindent {\bf Theorem.} {\sl Let $\partial$ be an sh-Lie algebra
structure on $\fra g$, that is, a coalgebra perturbation of the
differential $d$ on $\Sigm^{\mathrm c}[s \fra g]$. Then the given
contraction and the sh-Lie algebra structure $\partial$ on $\fra
g$ determine an sh-Lie algebra structure on $M$, that is, a
coalgebra perturbation $\ppartial$ of the coalgebra differential
$\dd^0$ on $\Sigm^{\mathrm c}[sM]$, a Lie algebra twisting cochain
\begin{equation*}
\tau \colon \Sigm_{\ppartial}^{\mathrm c}[sM] \longrightarrow
\mathcal L \Sigm_{\partial}^{\mathrm c}[s \fra g]
\end{equation*}
and, furthermore, a contraction
\begin{equation*}
   \left(\Sigm_{\ppartial}^{\mathrm c}[sM]
     \begin{CD}
      \null @>{\overline \tau}>> \null\\[-3.2ex]
      \null @<<{\Pi_{\partial}}< \null
     \end{CD}
    \mathcal C [\mathcal L \Sigm_{\partial}^{\mathrm c}[s \fra g]],
    H_{\partial}\right)
  \end{equation*}
of chain complexes which are natural in terms of the data.
The injection 
\[
\overline \tau
\colon \Sigm_{\ppartial}^{\mathrm c}[sM]\to \mathcal C[\mathcal L
\Sigm_{\partial}^{\mathrm c}[s \fra g]]
\]
is then a {\em
morphism\/} of coaugmented differential graded coalgebras.}

\smallskip\noindent
Together with the adjoint
$ \Sigm_{\partial}^{\mathrm c}[s \fra g]
\to \mathcal C[\mathcal L
\Sigm_{\partial}^{\mathrm c}[s \fra g]]$
of the universal Lie algebra twisting cochain of 
$\mathcal L
\Sigm_{\partial}^{\mathrm c}[s \fra g]$,
this yields an sh-equivalence between
$(M,\ppartial)$ and $(\fra g,\partial)$.
For the special case where $M$ and $\fra g$ are connected, 
we also construct an explicit extension of the retraction $\Pi_{\partial}$
to an sh-Lie map.}
\end{abstract}

\MSC 16E45, 16W30, 17B55, 17B56, 17B65, 17B70, 18G10, 55P62
\EndMSC

\section{Introduction}

{\em Higher homotopies\/} are nowadays playing a prominent role
in mathematics as well as in certain branches of theoretical
physics. Higher homotopies often arise as follows:
Suppose we are given a huge object, e.~g. a chain complex,
whose homology includes invariants of a certain geometric or 
algebraic situation.
When one tries to cut such a huge object to size
by passing to a smaller object, chain equivalent to the initial one,
typically higher homotopies, e.~g. {\em Massey products\/}, arise. 
Furthermore,
under  homotopy, strict
algebraic structures such as e.~g. 
the Jacobi identity of 
a differential graded Lie bracket
are not in general preserved,
and {\em higher homotopies\/} arise measuring e.~g. the failure
of the Jacobi identity in a coherent way.
Even for strict structures,
non-trivial higher homotopies may encapsulate
additional information; this is true, e.~g., for the Borromean rings:
A non-trivial Massey product detects the non-trivial linking of the 
three rings.
In physics such higher homotopies 
arise e.~g. as anomalies or
higher order correlation functions;
see e.~g. \cite{jimmurra} and the references there, in particular
to the seminal papers of J. Stasheff.

The ordinary perturbation lemma for chain complexes 
has become a standard tool to
handle higher homotopies in a constructive manner. 
In view of a celebrated result of Kontsevich's, 
sh-Lie (also known as $L_{\infty}$) algebras
have attracted much attention, and the issue of
compatibility of the perturbation lemma with a general sh-Lie 
 algebra structure
arises. The question whether
certain perturbation constructions preserve algebraic
structure actually shows up already when one tries to
construct e.~g. models for differential graded algebras.
In the literature, the {\em tensor trick\/} \cite{homotype},
\cite{huebkade}, cf. \cite{jimmurra} for more literature, 
was successfully exploited 
to explore perturbations of free differential graded algebras
and cofree differential graded coalgebras,
the basic reason for that success
being the fact that homotopies of 
morphisms of such algebras or coalgebras can then be handled concisely; 
this tensor trick may actually be viewed as an instance of
a labelled rooted trees construction \cite{tornike}.
However, for differential graded  cocommutative coalgebras
as well as for differential graded  commutative algebras, 
the tensor trick breaks down; 
indeed, as noted already in
\cite{schlstas}, the notion of {\em homotopy of morphisms of
cocommutative coalgebras is a subtle concept\/}.
The Cartan-Chevalley-Eilenberg coalgebra
(or classifying coalgebra) of a differential graded Lie algebra is a
differential graded cocommutative coalgebra;
more generally, an sh-Lie algebra is defined
in terms of a coalgebra perturbation on
a differential graded cocommutative coalgebra.
These objects actually arise
in deformation theory, see e.~g. \cite{berikash} and the literature there.
The purpose of the present paper is to 
offer ways to overcome 
the difficulties with the notion of homotopy
in the (co)commutative case
by establishing the {\em perturbation
lemma\/} for sh-Lie algebras. 
As a side remark we note that, in a different context, suitable 
homological perturbation theory (HPT)
constructions that are compatible
with other algebraic structure enabled us to carry out complete
numerical calculations in group cohomology \cite{cohomolo}--\cite{holo} 
which cannot be done by other
methods.

To explain this general perturbation
lemma at the present stage somewhat informally, let $R$ be a
commutative ring with 1 which contains the rational numbers
as a subring and let $(M
     \begin{CD}
      \null @>{\nabla}>> \null\\[-3.2ex]
      \null @<<{\pi}< \null
     \end{CD}
    \fra g, h)$
be a {\em contraction\/} of chain complexes over $R$. 
Differential graded Lie algebras defined over a ring more
general than a field arise in homotopy theory via Samelson brackets, 
cf. e.~g. \cite{comoonei},
in gauge theory, e.~g. as Lie algebras of gauge 
transformations---here the ground ring is the algebra of 
smooth functions on a smooth manifold
and hence manifestly contains the rationals as a 
subring---and in combinatorial group theory \cite{magkarso}.
These remarks justify, perhaps, building the theory
over rings more general than a field.
A version of the sh-Lie
algebra perturbation lemma is the following.

\medskip
\noindent {\bf Theorem.} {\sl Given an sh-Lie algebra structure on
$\fra g$, that is, a coalgebra perturbation of the differential
$d$ on $\Sigm^{\mathrm c}[s \fra g]$,
the chain complex $M$ acquires an sh-Lie algebra
structure that is natural in terms of the given contraction and
the sh-Lie algebra structure on $\fra g$, and the data determine
an sh-equivalence between $M$ and $\fra g$ relative to the sh-Lie
algebra structures that is natural in terms of the data. }
\medskip

The meaning of  sh-equivalence is this: Given the coalgebra perturbation
$\partial$ of the differential
$d$ on $\Sigm^{\mathrm c}[s \fra g]$,
the data determine in particular
a coalgebra perturbation $\ppartial$ of the coalgebra differential
$\dd^0$ on $\Sigm^{\mathrm c}[sM]$ and a Lie algebra twisting cochain
\begin{equation*}
\tau \colon \Sigm_{\ppartial}^{\mathrm c}[sM] \longrightarrow
\mathcal L \Sigm_{\partial}^{\mathrm c}[s \fra g].
\end{equation*}
The injection 
$
\overline \tau
\colon \Sigm_{\ppartial}^{\mathrm c}[sM]\to \mathcal C[\mathcal L
\Sigm_{\partial}^{\mathrm c}[s \fra g]]
$
is then a morphism of coaugmented differential graded coalgebras
inducing an isomorphism on homology.
Together with the adjoint
$ \Sigm_{\partial}^{\mathrm c}[s \fra g]
\to \mathcal C[\mathcal L
\Sigm_{\partial}^{\mathrm c}[s \fra g]]$
of the universal Lie algebra twisting cochain of 
$\mathcal L
\Sigm_{\partial}^{\mathrm c}[s \fra g]$,
this yields an sh-equivalence between
$(M,\ppartial)$ and $(\fra g,\partial)$.

A special case of the theorem is the Lie algebra perturbation lemma
established in a predecessor of this paper \cite{pertlie}. 
Exploiting 
a suitable version of the loop Lie algebra
relative to a coaugmented differential graded cocommutative coalgebra,
see Section 2 below for details,
we will reduce the present general case to the special case in
\cite{pertlie}. 
We conjecture that the theory we develop in this paper
has applications to foliation theory
and to the integration problem of sh-Lie algebras. 
The main result of the  present paper includes a very general
solution of the {\em master equation\/} or, equivalently,
{\em Maurer-Cartan equation\/}.
More comments about the relevance and history of
the master equation can be found in 
\cite{pertlie}, \cite{jimmurra}, and \cite{huebstas}.

I am much indebted to Jim Stasheff for having prodded me on
various occasions to pin down the general perturbation lemma for
sh-Lie algebras, to M. Duflo for discussion about
the PBW theorem,  and to the referee for 
a number of comments which helped improve the exposition.

\section{The sh-Lie algebra perturbation lemma}
\label{outline}

The ground ring, written as $R$, is a commutative ring with $1$ 
which contains the rationals as a subring. 
We will take {\em chain complex\/} to mean {\em
differential graded\/} $R$-{\em module\/}. A chain complex will
not necessarily be concentrated in non-negative or non-positive
degrees. The differential of a chain complex will always be
supposed to be of degree $-1$. Write $s$ for the {\em
suspension\/} operator as usual and, accordingly, $s^{-1}$ for the
{\em desuspension\/} operator. Thus, given the chain complex $X$,
$(sX)_j = X_{j-1}$, etc., and the differential ${d\colon sX \to
sX}$ on the suspended object $sX$ is defined in the standard
manner so that $ds+sd=0$.

For a filtered chain complex $X$, a {\em perturbation\/} of the
differential $d$ of $X$ is a (homogeneous) morphism $\partial$ of
the same degree as $d$ such that $\partial$ lowers filtration and
$(d +
\partial)^2 = 0$ or, equivalently,
\begin{equation}
[d,\partial] + \partial \partial = 0.
\end{equation}
Thus, when $\partial$ is a perturbation on $X$, the sum $d +
\partial$, referred to as the {\em perturbed differential\/},
endows $X$ with a new differential. When $X$ has a graded
coalgebra structure such that $(X,d)$ is a differential graded
coalgebra, and when the {\em perturbed differential\/} $d +
\partial$ is compatible with the graded coalgebra structure, we
refer to $\partial$ as a {\em coalgebra perturbation\/}; the
notion of {\em algebra perturbation\/} is defined accordingly.
Given a differential graded coalgebra $C$ and a coalgebra
perturbation $\partial$ of the differential $d$ on $C$, we will
occasionally denote the new differential graded coalgebra by
$C_{\partial}$. Thus the differential of the latter is given by
the sum $d +\partial$.

The following notion goes back to \cite{eilmactw}: A {\em
contraction\/}
\begin{equation}
   (N
     \begin{CD}
      \null @>{\nabla}>> \null\\[-3.2ex]
      \null @<<{\pi}< \null
     \end{CD}
    M, h) \label{co}
  \end{equation}
of chain complexes  consists of
chain complexes $N$ and $M$,
chain maps $\pi\colon N \to M$ and $\nabla \colon M \to
N$,
and a morphism $h\colon N \to N$ of the underlying graded
modules of degree 1;
these data are required to satisfy
\begin{align}
 \pi \nabla &= \mathrm{Id},
\label{co0}
\\
Dh &= \mathrm{Id} -\nabla \pi, \label{co1}
\\
\pi h &= 0, \quad h \nabla = 0,\quad hh = 0. \label{side}
\end{align}
The requirements \eqref{side} are referred to as {\em annihilation
properties\/} or {\em side conditions\/}.

\begin{Remark} \label{rem1} It is well known that the side 
conditions \eqref{side}
can always be achieved. This fact relies on the standard observation that
a chain complex is contractible if and only
if it is isomorphic to a cone, cf. {\rm \cite{husmosta} (IV.1.5)}.
Under the present circumstances,
given data of the kind \eqref{co} such
that \eqref{co0} and \eqref{co1} hold but not necessarily the side
conditions \eqref{side}, the operator
\[
\widetilde h=(\mathrm{Id} -\nabla \pi)h(\mathrm{Id} -\nabla
\pi)d(\mathrm{Id} -\nabla \pi)h(\mathrm{Id} -\nabla \pi)
\]
satisfies the requirements \eqref{co1} and \eqref{side}, with
$\widetilde h$ instead of $h$;
when $h$ already satisfies \eqref{side}, $\widetilde h$ coincides
with $h$.
\end{Remark}

Let $C$ be a {\em coaugmented\/} differential graded coalgebra
with coaugmentation map $\eta \colon R \to C$ and {\em
coaugmentation\/} coideal $JC = \mathrm{coker}(\eta)$, the
diagonal map being written as $\Delta \colon C \to C \otimes C$. 
Recall that the counit $\varepsilon \colon C \to R$ and the
coaugmentation map determine a direct sum decomposition $C = R
\oplus JC$ and that the {\em coaugmentation\/} filtration $\{\mathrm
F_nC\}_{n \geq 0}$ is given by
\[\mathrm F_nC = \mathrm{ker}(C \longrightarrow (JC)^{\otimes
(n+1)})\  (n \geq 0)
\]
where the unlabelled arrow is induced by some iterate of the
diagonal $\Delta$ of $C$. This filtration is well known to turn
$C$ into a {\em filtered\/} coaugmented differential graded
coalgebra; thus, in particular, $\mathrm F_0C = R$. We recall that
$C$ is said to be {\em cocomplete\/} when $C=\cup \mathrm F_nC$.

Let $C$ be a coaugmented differential graded coalgebra  and $A$ an
augmented differential graded algebra, the 
multiplication map of $A$ being written as 
$\mu \colon A \otimes A \to A$ and the augmentation map as
$\varepsilon \colon A \to R$. Recall that, given
homogeneous morphisms $a,b \colon C \to A$, their {\em cup
product\/} $a\cup b$ is the composite
\begin{equation}
\begin{CD} C @>{\Delta}>> C\otimes C @>{a\otimes b}>>A
\otimes´A @> {\mu}>> A
\end{CD}
\label{cupproduct}
\end{equation}
where $\mu$ refers to the multiplication map of $A$. The cup product
$\cup$ is well known to turn $\mathrm{Hom}(C,A)$ into an augmented
differential graded algebra, the differential being the ordinary
Hom-differential. Recall also that an {\em ordinary twisting
cochain\/}
$\tau \colon C \longrightarrow A$
is a homogeneous morphism of the underlying graded $R$-modules of
degree $-1$ satisfying the identity
\begin{equation}
D\tau = \tau \cup \tau \label{tc6}
\end{equation}
and the requirements $\tau \eta=0$ and $\varepsilon \tau= 0$.

Given
two graded objects $U$ and $V$,  we denote the (graded)
interchange map by
$T \colon U \otimes V \longrightarrow V \otimes U$.
Recall that a graded coalgebra $C$ is graded cocommutative
when its diagonal map $\Delta$ satisfies the condition $T\Delta = \Delta$.

Let $\fra g$ be (at first) a
chain complex, the differential being written as $d \colon \fra g
\to \fra g$, and let
\begin{equation}
   (M
     \begin{CD}
      \null @>{\nabla}>> \null\\[-3.2ex]
      \null @<<{\pi}< \null
     \end{CD}
    \fra g, h)
\label{cont1}
  \end{equation}
be a {\em contraction\/} of chain complexes. 
Consider 
the {\em cofree\/} coaugmented
differential graded {\em cocommutative\/} coalgebra
(differential graded {\em symmetric\/} coalgebra) 
${\Sigm^{\mathrm c} = \Sigm^{\mathrm c}[sM]}$ on
the suspension $sM$ of $M$, the existence of 
of that coalgebra being guaranteed by the hypothesis that
the ground ring $R$ contain the rational numbers as a subring.
Further, let
$\dd^0$ denote the coalgebra differential on $\Sigm^{\mathrm c} =
\Sigm^{\mathrm c}[sM]$ induced by the differential on $M$. For $b
\geq 0$, we will henceforth denote the homogeneous
(tensor) degree $b$
component of $\Sigm^{\mathrm c}[sM]$ by $\Sigm_b^{\mathrm c}$;
thus, as a chain complex, $\mathrm F_b\Sigm^{\mathrm c} = R \oplus
\Sigm_1^{\mathrm c} \oplus \dots \oplus \Sigm_b^{\mathrm c}$.
Likewise, as a chain complex, $\Sigm^{\mathrm c} =
\oplus_{j=0}^{\infty} \Sigm_j^{\mathrm c}$. We denote by
$
\tau_{M}\colon \Sigm^{\mathrm c} \longrightarrow M
$
the composite of the canonical projection $\mathrm{proj}\colon
\Sigm^{\mathrm c} \to sM$ from $\Sigm^{\mathrm c} = \Sigm^{\mathrm
c}[sM]$ to its homogeneous degree 1 constituent $sM$ with
the desuspension map $s^{-1}$ from $sM$ to $M$. 
In particular,
$\tau_{\fra g}\colon \Sigm^{\mathrm c}[s\fra g]\longrightarrow \fra
g
$
refers to the composite of the canonical projection to $\Sigm_1^{\mathrm
c}[s\fra g] =s\fra g$ with the desuspension map.

Given a homogeneous element $x$ of a graded module, we will
denote its degree by $|x|$. Given two chain complexes $X$ and $Y$, recall
that $\mathrm{Hom}(X,Y)$ inherits the structure of a chain complex
by the operator $D$ defined by
$D \phi = d \phi -(-1)^{|\phi|} \phi d$
where $\phi$ is a homogeneous morphism of $R$-modules from $X$ to $Y$.

Let now  $C$ be a coaugmented differential
graded cocommutative coalgebra
and $\fra h$ a differential graded Lie algebra, the graded bracket
being written as $[\, \cdot \, , \, \cdot \, ]$.
Given homogeneous morphisms $a,b
\colon C \to \fra h$, with a slight abuse of the bracket notation
$[\, \cdot \, , \, \cdot \, ]$, their {\em cup bracket\/} $[a, b]$
is given by the composite
\begin{equation}
\begin{CD} C @>{\Delta}>> C\otimes C @>{a\otimes b}>>\fra h
\otimes\fra h @> {[\cdot,\cdot]}>> \fra h.
\end{CD}
\label{cupbracket}
\end{equation}
This bracket turns $\mathrm{Hom}(C,\mathfrak h)$ 
into a differential graded Lie algebra.

In particular, take $C$ to be
the differential graded symmetric
coalgebra 
$\Sigm^{\mathrm c}[s\fra h]$ and
define the coderivation
\begin{equation}
\partial\colon\Sigm^{\mathrm
c}[s\fra h] \longrightarrow \Sigm^{\mathrm c}[s\fra h] \label{cod}
\end{equation}
on $\Sigm^{\mathrm c}[s\fra h]$  by the requirement
\begin{equation}
\tau_{\fra h} \partial = \frac 12 [\tau_{\fra h}, \tau_{\fra
h}]\colon \Sigm_2^{\mathrm c}[s\fra h] \to \fra h. \label{proc333}
\end{equation}
Plainly $D\partial\  (=d\partial + \partial d) = 0$ since the Lie
algebra structure on $\fra h$ is supposed to be compatible with
the differential $d$ on $\fra h$. Moreover, the 
requirement that the bracket $[\,
\cdot \, , \, \cdot \, ]$ on $\fra h$ satisfy the graded Jacobi
identity is equivalent to the requirement that $\partial\partial$ vanish, 
that is, to $\partial$ being a coalgebra perturbation of the differential $d^0$
on $\Sigm^{\mathrm c}[s\fra h]$, cf. \cite{pertlie} and
\cite{huebstas}. 
The  Lie algebra perturbation lemma
(Theorem 2.1 in \cite{pertlie} and reproduced below as Lemma
\ref{lem11}) and the sh-Lie algebra perturbation lemma (Theorem
\ref{lem13} below) both generalize this observation. 
Under the present circumstances, $\fra h$ being an ordinary differential
graded Lie algebra, the resulting
differential graded coalgebra $\Sigm_{\partial}^{\mathrm c}[s\fra
h]$ is precisely the standard {\scs Cartan-Chevalley-Eilenberg\/}
(CCE-) or {\em classifying\/} coalgebra for $\fra h$ and,
following \cite{quilltwo} (p.~291), we denote this coalgebra by $\mathcal
C[\fra h]$.

As before, let $C$ be a coaugmented differential graded cocommutative coalgebra.
A {\em Lie algebra twisting cochain\/} $t \colon C \to \fra h$ is
a homogeneous morphism of degree $-1$ whose composite with the
coaugmentation  map is zero and which
 satisfies the equation
\begin{equation}
Dt = \frac 12 [t,t] , \label{master}
\end{equation}
cf. \cite{huebstas}, \cite{moorefiv} and \cite{quilltwo}.
The equation \eqref{master}
is a version of the {\em master equation\/}, cf. \cite{huebstas}
and the literature there. In particular, relative to the graded
Lie bracket $\ppppartial$ on $\fra h$, the morphism  
$\tau_{\fra h}\colon \mathcal C[\fra h] \to \fra
h$ is a Lie algebra twisting cochain, the {\scs
Cartan-Chevalley-Eilenberg\/} (CCE-) or {\em universal\/} Lie
algebra twisting cochain for $\fra h$.  Likewise, when $M$ is
viewed as an {\em abelian\/} differential graded Lie algebra,
$\Sigm^{\mathrm c}=\Sigm^{\mathrm c}[sM]$ may be viewed as the
{\scs CCE-\/} or {\em classifying\/} coalgebra $\mathcal C [M]$
for $M$, and $\tau_{M} \colon \Sigm^{\mathrm c} \to M$ is then the
universal differential graded Lie algebra twisting cochain for $M$
as well.

For intelligibility, we will now recall the main result of
\cite{pertlie}, spelled out there as Theorem 2.1.

\begin{Lemma}[Lie algebra perturbation lemma]
\label{lem11} Suppose that $\fra g$ carries a differential graded
Lie algebra structure. Then the contraction {\rm {\eqref{cont1}}}
and the graded Lie algebra structure on $\fra g$ determine an
sh-Lie algebra structure on $M$, that is, a coalgebra perturbation
$\ppartial$ of the coalgebra differential $\dd^0$ on
$\Sigm^{\mathrm c}[sM]$, a Lie algebra twisting cochain
\begin{equation}
\tau \colon \Sigm_{\ppartial}^{\mathrm c}[sM] \longrightarrow \fra
g \label{tc1}
\end{equation}
and, furthermore, a contraction
\begin{equation}
   \left(\Sigm_{\ppartial}^{\mathrm c}[sM]
     \begin{CD}
      \null @>{\overline \tau}>> \null\\[-3.2ex]
      \null @<<{\Pi}< \null
     \end{CD}
    \mathcal C[\fra g], H\right)
\label{cont5}
  \end{equation}
of chain complexes which are natural in terms of the data so that
\begin{align}
\pi \tau&=\tau_M\colon \Sigm^{\mathrm c}[sM] \longrightarrow M,
\label{twist33}
\\
h \tau &= 0.
\label{twist44}
\end{align}
The injection $\overline \tau \colon \Sigm_{\ppartial}^{\mathrm c}[sM]\to
\mathcal C[\fra g]$ 
is then a {\em morphism\/} of
coaugmented differential graded coalgebras.
\end{Lemma}

In the statement of Lemma \ref{lem11}, 
the adjoint $\overline \tau$ of \eqref{tc1} is plainly an
sh-equivalence in the sense that it induces an {\em isomorphism on
homology\/}, including the brackets of all order that are induced on
homology, $M$ being endowed with the sh-Lie algebra structure given by
$\ppartial$. In Section \ref{invert} below we shall explain how
$\tau$ yields actually an sh-equivalence between $\fra g$ and $M$
in a certain stronger sense when $M$ and $\fra g$ are connected.

We will now consider the more general case where $\fra g$ is endowed
with merely an sh-Lie algebra structure. 
To this end, we will 
denote by $\Sigm$ the graded symmetric algebra functor in the
category of $R$-modules. As before, let $\fra h$ be a differential graded Lie
algebra and, $\fra h$ being viewed as a chain complex, let
$\Sigm[\fra h]$ be the differential graded symmetric algebra on
$\fra h$. Since the ground $R$ is supposed to contain the rational
numbers as a subring,
the diagonal map $\mathfrak h
\to \mathfrak h \oplus \mathfrak h$ of $\mathfrak h$ 
induces a diagonal map $\Delta\colon \Sigm[\fra h]\to
\Sigm[\fra h] \otimes \Sigm[\fra h]$ that turns 
$\Sigm[\fra h]$
into a differential graded cocommutative Hopf
algebra; furthermore, the obvious filtration
then turns 
$\Sigm[\fra h]$
into a filtered differential graded cocommutative Hopf
algebra.
Consider the universal differential graded algebra $\mathrm
U[\fra h]$ associated with $\fra h$ and let $j \colon \fra h \to
\mathrm U[\fra h]$ denote the canonical morphism of differential
graded Lie algebras; it is well known that, 
just as for the symmetric algebra on $\fra h$, 
via the appropriate
universal property, the diagonal map of $\mathfrak h$
induces a
diagonal map $\Delta\colon\mathrm U[\fra h] \to \mathrm U[\fra
h]\otimes \mathrm U[\fra h]$ turning $\mathrm U[\fra h]$ into a
differential graded cocommutative Hopf algebra 
which, relative to the
ordinary Poincar\'e-Birkhoff-Witt filtration
$
R=\mathrm U_0\subseteq \mathrm U_1 \subseteq \ldots \subseteq
\mathrm U_{\ell} \subseteq \ldots
$
is actually a filtered
differential graded cocommutative Hopf algebra. We denote the
associated graded object by $\mathrm E^0\mathrm U[\fra h]$; this
is a differential graded commutative and cocommutative
Hopf algebra endowed with
a canonical  morphism
$
\Sigm[\fra h] \longrightarrow \mathrm E^0\mathrm U[\fra h]
$
of differential graded Hopf algebras. 

\begin{Proposition}
\label{PBW}
The classical Poincar\'e symmetrization map
\[
e\colon \Sigm[\fra h] \longrightarrow \mathrm U[\fra h],\quad
e(x_1\ldots x_n)=\frac 1{n!}\sum _{\sigma\in S_n} \pm j(x_{\sigma 1})
\ldots j(x_{\sigma n})
\]
in the category of $R$-modules
is a functorial isomorphism of filtered differential graded ($R$-)coalgebras
which induces an isomorphism
\[
e\colon \Sigm[\fra h] \longrightarrow \mathrm E^0\mathrm U[\fra h]
\]
of differential graded Hopf algebras.
Consequently, the differential graded ($R$-)algebra $\mathrm U[\fra h]$
being viewed as a differential graded  ($R$-)Lie algebra via 
the commutator bracket
as usual,
the canonical morphism
 $j \colon \fra h \to
\mathrm U[\fra h]$  of differential
graded $R$-Lie algebras is injective, and
the universal algebra $\mathrm U[\fra h]$ is enveloping. \qed
\end{Proposition}

This proposition makes precise the idea
that $\mathrm U[\fra h]$, viewed as a differential graded
Hopf algebra, is a {\em
perturbation\/} of  $\Sigm[\fra h]$, viewed as a differential graded
Hopf algebra, the coalgebra structure being unperturbed.

Let $Y$ be a chain complex, and let $\mathrm T [Y]$ be the
differential graded tensor algebra on $Y$. The {\em shuffle\/}
diagonal map is well known to turn $\mathrm T[Y]$ into a
differential graded cocommutative Hopf algebra and, $\mathrm T
[Y]$ being viewed as a differential graded Lie algebra via the
commutator bracket, the {\em free\/} (differential graded) Lie
algebra $\mathrm L[Y]$ on $Y$ is the differential graded Lie
subalgebra of $\mathrm T [Y]$ generated by $Y$.  Further, the
canonical morphism
of augmented differential graded algebras 
from $\mathrm U[\mathrm L[Y]]$ to  $\mathrm T[Y]$
is an
isomorphism, cf. e.~g. \cite{comoonei} (Proposition 2.10). 
This explains the
differential graded Hopf algebra structure on
$\mathrm U[\mathrm L[Y]]$ in the particular case of the differential
graded Lie algebra $\mathrm L[Y]$.

The submodule $\mathrm{Prim}[Y]$ of {\em primitive elements\/} in
the Hopf algebra $\mathrm T[Y]$ is well known to be a differential
graded Lie subalgebra of $\mathrm T[Y]$ and, since $Y$ is
manifestly contained in $\mathrm{Prim}[Y]$, the {\em free\/}
(differential graded) Lie algebra $\mathrm L[Y]$ on $Y$ is plainly
a differential graded Lie subalgebra of $\mathrm{Prim}[Y]$. In
view of a classical theorem of K.~O. Friedrichs', over a field of
characteristic zero, the two coincide and, more generally, the two
coincide whenever the ground ring $R$ is an integral domain of
characteristic zero and $Y$  a free graded $R$-module, cf.
\cite{comoonei} (Proposition 2.8).

Let $C$ be a coaugmented differential graded coalgebra.
By construction, the loop algebra $\Omega C$ is the {\em
perturbed\/} tensor algebra $\mathrm T_{\Delta} [s^{-1}JC]$ on
$s^{-1}JC$, the {\em algebra perturbation\/} $\partial_{\Delta}$
on $\mathrm T[s^{-1}JC]$ being induced by $\Delta$.
Suppose,
in addition, that $C$ is cocommutative. Then $\Omega C$ acquires a
differential graded Hopf algebra structure. Moreover, since the
diagonal map $\Delta$ is a morphism of differential graded
coalgebras, the induced morphism 
$J\Delta\colon JC \to JC \otimes JC$
is compatible with the structure, whence the
algebra perturbation $\partial_{\Delta}$ descends to a {\em Lie
algebra perturbation\/} on
\[
\mathrm{Prim}[s^{-1}JC] = \mathrm {ker}(J\Delta)
\]
which we still denote by $\partial_{\Delta}$, and we denote the
resulting differential graded Lie algebra by
$\mathrm{Prim}_{\Delta}[s^{-1}JC]$. Over a field of characteristic
zero, this is the {\em loop\/} Lie algebra on $C$, a familiar object,
and the loop Lie algebra then coincides with the free Lie algebra.
In general, it is not clear whether
the obvious injection of
the
free differential graded Lie algebra $\mathrm{L}[s^{-1}JC]$ into
$\mathrm{Prim}[s^{-1}JC]$
is onto.

\begin{Lemma} The coaugmented differential graded
coalgebra $C$ being assumed to be graded cocommutative, 
the values of the Lie
algebra perturbation $\partial=\partial_{\Delta}$, restricted to
$\mathrm{L}[s^{-1}JC]$, lie in $\mathrm{L}[s^{-1}JC]$.
\end{Lemma}

\begin{proof} Write $Y=s^{-1}JC$, so that
$\mathrm{L}[s^{-1}JC]=\mathrm{L}[Y]\subseteq \Omega C$, and so
that the augmented graded algebra which underlies $\Omega C$
coincides with the tensor algebra $\mathrm T[Y]$.
We will use the notation
$[\,\cdot \, , \, \cdot \,]$ for the graded commutator in the
graded tensor algebra $\mathrm T[s^{-1}JC]$.
The values of
the morphism
\[
\partial - T\partial \colon Y \longrightarrow Y \otimes Y
\]
lie in the submodule $[Y,Y]\subseteq Y \otimes Y$ spanned by the
commutators of elements from $Y$.  The {\em algebra
perturbation\/} $\partial_{\Delta}$ on $\mathrm T[s^{-1}JC]$ is
induced by the morphism $J\Delta$ coming from the diagonal map
$\Delta$ of $C$. Since the latter is cocommutative, $- T\partial$
coincides with $\partial$ whence the values of $2 \partial$,
restricted to $Y$, lie in $\mathrm{L}[Y]$. 
Consequently the values of the
perturbation $\partial$, restricted to $Y$, lie in
$\mathrm{L}[Y]$.
\end{proof}

Let $C$ be a coaugmented differential graded cocommutative coalgebra.
We will  use the notation $\mathcal L C$ for
$\mathrm{L}[s^{-1}JC]$,
endowed with the perturbed differential $d+\partial_{\Delta}$ 
and we will refer to $\mathcal L C$ as the {\em loop\/} Lie algebra over
$C$. 
The desuspension map induces a Lie algebra twisting
cochain
\[
t_{\mathcal L} \colon C \longrightarrow \mathcal LC,
\]
the {\em universal Lie algebra twisting cochain\/} for the loop
Lie algebra. See \cite{moorefiv} and \cite{quilltwo} for the case
where the ground ring is a field of characteristic zero.
Whether or not the ground ring is a field of characteristic zero,
the canonical morphism
\begin{equation}
\mathrm U[\mathcal L C] \longrightarrow \Omega C
\label{canalg}
\end{equation}
of augmented differential graded algebras is an isomorphism,
and the adjoint
\begin{equation*}
\Omega C \longrightarrow \mathrm U[\mathcal L C]
\end{equation*}
of the composite of $t_{\mathcal L}$ with the canonical morphism
$\mathcal LC \to \mathrm U[\mathcal LC]$ yields the inverse for
\eqref{canalg} in the category of augmented differential graded
algebras.

With $C=\Sigm^{\mathrm c}[s\fra g]$, the
isomorphism \eqref{canalg} then takes the form
\begin{equation}
\mathrm U[\mathcal L \Sigm^{\mathrm c}[s\fra g]] \longrightarrow
\Omega \Sigm^{\mathrm c}[s\fra g] . \label{canalg2}
\end{equation}
An sh-{\em Lie algebra structure\/} or $L_{\infty}$-{\it
structure\/} on the chain complex $\fra g$ is a {\em coalgebra
perturbation\/} $\partial$ of the differential $d$ on the cofree
coaugmented differential graded cocommutative coalgebra
$\Sigm^{\mathrm c}[s \fra g]$ on $s \fra g$, cf. \cite{huebstas}
(Def.~2.6). Given such an sh-Lie algebra structure $\partial$ on
$\fra g$, with $C=\Sigm_{\partial}^{\mathrm c}[s\fra g]$, 
\eqref{canalg} yields the isomorphism
\begin{equation}
\mathrm U[\mathcal L \Sigm_{\partial}^{\mathrm c}[s\fra g]]
\longrightarrow \Omega \Sigm_{\partial}^{\mathrm c}[s\fra g] .
\label{canalg3}
\end{equation}
In particular, via the coderivation \eqref{cod}, an ordinary
graded Lie algebra structure $\ppppartial$ 
determines an sh-Lie algebra structure $\partial$ and, in this
case, $\Sigm_{\partial}^{\mathrm c}[s\fra g]$ amounts to
the {\sc CCE\/}-coalgebra $\mathcal C[\fra g]$ for $(\fra g,\ppppartial)$.
Given two sh-Lie
algebras $(\fra g_1,
\partial_1)$ and $(\fra g_2,\partial_2)$, an {\em sh-morphism\/}
or {\em sh-Lie map\/} from $(\fra g_1, \partial_1)$ to $(\fra g_2,
\partial_2)$ is a morphism ${ \Sigm^{\mathrm c}_{\partial_1}[s\fra
g_1] \to \Sigm^{\mathrm c}_{\partial_2}[s\fra g_2] }$ of
coaugmented differential graded coalgebras \cite{huebstas};
 we then define a {\em generalized sh-morphism\/} or
{\em generalized sh-Lie map\/} from $(\fra g_1, \partial_1)$ to
$(\fra g_2, \partial_2)$ to be  a Lie algebra twisting cochain ${
\Sigm^{\mathrm c}_{\partial_1}[s\fra g_1] \to \mathcal L
\Sigm^{\mathrm c}_{\partial_2}[s\fra g_2] }$.

\begin{Theorem}[Sh-Lie algebra perturbation lemma]
\label{lem13} Let $\fra g$ be a chain complex and
let $\partial$ be an sh-Lie algebra structure on $\fra g$, that
is, a coalgebra perturbation of the differential $d$ on
$\Sigm^{\mathrm c}[s \fra g]$. Then the contraction {\rm
{\eqref{cont1}}} and the sh-Lie algebra structure $\partial$ on
$\fra g$ determine an sh-Lie algebra structure on $M$, that is, a
coalgebra perturbation $\ppartial$ of the coalgebra differential
$\dd^0$ on $\Sigm^{\mathrm c}[sM]$, a Lie algebra twisting cochain
\begin{equation}
\tau \colon \Sigm_{\ppartial}^{\mathrm c}[sM] \longrightarrow
\mathcal L \Sigm_{\partial}^{\mathrm c}[s \fra g] \label{tc2}
\end{equation} and, finally, a contraction
\begin{equation}
   \left(\Sigm_{\ppartial}^{\mathrm c}[sM]
     \begin{CD}
      \null @>{\overline \tau}>> \null\\[-3.2ex]
      \null @<<{\Pi_{\partial}}< \null
     \end{CD}
    \mathcal C [\mathcal L \Sigm_{\partial}^{\mathrm c}[s \fra g]],
    H_{\partial}\right) \label{cont555}
  \end{equation}
of chain complexes, and
{\rm \eqref{tc2}}
and
{\rm \eqref{cont555}}
are natural in terms of the data.
The injection 
\[\overline \tau
\colon \Sigm_{\ppartial}^{\mathrm c}[sM]\to \mathcal C[\mathcal L
\Sigm_{\partial}^{\mathrm c}[s \fra g]] 
\]
is then a {\em
morphism\/} of coaugmented differential graded coalgebras.
\end{Theorem}

Under the circumstances of Theorem \ref{lem13}, the twisting
cochain \eqref{tc2} is a generalized morphism of sh-Lie algebras
from $(M,\pppartial)$ to $(\fra g, \partial)$, and the adjoint
$\overline \tau$ of \eqref{tc2} is plainly an sh-equivalence in
the sense that it induces an isomorphism on homology, including
the brackets of all order that are induced on homology. In Section
\ref{invert} below, we shall sketch an extension of the
contraction \eqref{cont555} to an sh-equivalence, in a stronger
sense, between these two sh-Lie algebras for the special case where 
$M$ and $\fra g$ are connected.

\section{Proof of the sh-Lie algebra perturbation lemma}
\label{proof}

Until further notice we will view $\fra g$ merely as a chain
complex or, equivalently, as an {\em abelian\/} differential
graded Lie algebra. The desuspension map induces the 
ordinary twisting cochain
\[
\tau^{\Sigm^{\mathrm c}}\colon \Sigm^{\mathrm c}[s \fra g]
\longrightarrow \Sigm [\fra g],
\]
and the adjoint $\pi_{\Sigm}\colon \Omega\Sigm^{\mathrm c}[s \fra
g] \to \Sigm[\fra g]$ thereof is a surjective morphism of
augmented differential graded algebras.

We will denote the {\em reduced\/} bar construction functor by
$\mathrm B$ (defined on
the category of augmented differential graded algebras).

\begin{Lemma}\label{LLL}
The projection $\pi_{\Sigm}$ extends to a contraction
\begin{equation}
   \left(\Sigm[\fra g]
     \begin{CD}
      \null @>{\nabla_{\Sigm}}>> \null\\[-3.2ex]
      \null @<<{\pi_{\Sigm}}< \null
     \end{CD}
    \Omega\Sigm^{\mathrm c}[s \fra g], h_{\Sigm}\right)
\label{cont51}
  \end{equation}
of chain complexes that is natural in terms of the data.
\end{Lemma}

In this lemma, nothing is claimed as far as compatibility of
$\nabla_{\Sigm}$ and $h_{\Sigm}$ with the algebra structures is concerned.

\begin{proof} Consider the ordinary loop algebra contraction
\begin{equation}
   \left(\Sigm[\fra g]
     \begin{CD}
      \null @>{\nabla^{\Omega}}>> \null\\[-3.2ex]
      \null @<<{\pi^{\Omega}}< \null
     \end{CD}
    \Omega\mathrm B\Sigm [\fra g], h^{\Omega}\right)
\label{cont511}
  \end{equation}
for $\Sigm[\fra g]$,  cf. \cite{husmosta}, \cite{munkholm} (2.14)
(p.~17). Here the projection $\pi^{\Omega}$ is the  adjoint of the
universal bar construction
twisting cochain $ \mathrm B\Sigm [\fra g]\to \Sigm[\fra
g]$ and is therefore a morphism of augmented differential graded
algebras. The adjoint
\begin{equation}
\nabla_{\Sigm^{\mathrm c}}=\overline{\tau^{\Sigm^{\mathrm
c}}}\colon\Sigm^{\mathrm c}[s \fra g] \longrightarrow \mathrm
B\Sigm [\fra g] \label{adj1}
\end{equation}
of the twisting cochain $\tau^{\Sigm^{\mathrm c}}$ is the standard
coalgebra injection of $\Sigm^{\mathrm c}[s \fra g]$ into $\mathrm B\Sigm
[\fra g]$, and a familiar construction extends \eqref{adj1} to a
contraction
\begin{equation}
   \left(\Sigm^{\mathrm c}[s\fra g]
     \begin{CD}
      \null @>{\nabla_{\Sigm^{\mathrm c}}}>> \null\\[-3.2ex]
      \null @<<{\pi_{\Sigm^{\mathrm c}}}< \null
     \end{CD}
    \mathrm B\Sigm[\fra g], h_{\Sigm^{\mathrm c}}\right)
\label{cont5111}
  \end{equation}
which is natural in terms of the data. Similarly, the induced
morphism
\begin{equation}
\Omega\nabla_{\Sigm^{\mathrm
c}}=\Omega\overline{\tau^{\Sigm^{\mathrm
c}}}\colon\Omega\Sigm^{\mathrm c}[s \fra g] \longrightarrow
\Omega\mathrm B\Sigm [\fra g] \label{adj2}
\end{equation}
of differential graded algebras extends to a contraction
\begin{equation}
   \left(\Omega\Sigm^{\mathrm c}[s\fra g]
     \begin{CD}
      \null @>{\Omega\nabla_{\Sigm^{\mathrm c}}}>> \null\\[-3.2ex]
      \null @<<{\pi_{\Omega\Sigm^{\mathrm c}}}< \null
     \end{CD}
    \Omega\mathrm B\Sigm[\fra g], h_{\Omega\Sigm^{\mathrm c}}\right)
\label{cont51111}
  \end{equation}
which is natural in terms of the data, and $\pi_{\Sigm} =
\pi^{\Omega}\circ \Omega\nabla_{\Sigm^{\mathrm c}}$. Let
\[
\nabla_{\Sigm}= \pi_{\Omega\Sigm^{\mathrm c}}
\circ\nabla^{\Omega},\ \widetilde h = \pi_{\Omega\Sigm^{\mathrm
c}}\circ h_{\Omega\Sigm^{\mathrm c}} \circ
\Omega\nabla_{\Sigm^{\mathrm c}},\ h_{\Sigm}=\widetilde h \circ d
\circ\widetilde h.
\]
This yields data of the kind \eqref{cont51}.
In view of Remark \ref{rem1} above, these data
constitute a
contraction of chain complexes that is natural in terms of the
data. \end{proof}

The chain complex $\fra g$ still being viewed as an
abelian differential graded Lie algebra, consider the loop Lie algebra
$\mathcal L= \mathcal L \Sigm^{\mathrm c}[s\fra g]$ on
$\Sigm^{\mathrm c}[s\fra g]$. Let $\nabla_{\mathcal
L}\colon \fra g \to \mathcal L\Sigm^{\mathrm c}[s \fra g]$ be the
canonical injection of chain complexes and, likewise, $\fra g$
still being viewed as an abelian differential graded Lie algebra,
let $\pi_{\mathcal L}\colon \mathcal L\Sigm^{\mathrm c}[s \fra g]
\to \fra g$ be the familiar adjoint of the corresponding universal
Lie algebra twisting cochain $\Sigm^{\mathrm c}[s \fra g]\to \fra
g$; this morphism $\pi_{\mathcal L}$ is plainly  a surjective
morphism of differential graded Lie algebras. It admits the
following elementary description: The canonical projection
$s^{-1}J\Sigm^{\mathrm c}[s \fra g]\to \fra g$ induces a
surjective morphism $\mathrm L[s^{-1}J\Sigm^{\mathrm c}[s \fra
g]]\to \mathrm L[\fra g]$ of differential graded Lie algebras, the
canonical projection $\mathrm L[\fra g]\to \fra g$ is simply the
abelianization map (of differential graded Lie algebras), and the
composite
\begin{equation}
\mathrm L[s^{-1}J\Sigm^{\mathrm c}[s \fra g]]\longrightarrow\fra g
\label{compo}
\end{equation}
of the two yields the morphism $\pi_{\mathcal L}$ of differential
graded Lie algebras, manifestly surjective, $\fra g$ being viewed abelian.

For intelligibility, we explain the details: Write $\mathrm L=
\mathrm L[s^{-1}J\Sigm^{\mathrm c}[s \fra g]]$ and let $\widetilde
{\mathrm L}$ denote the kernel of \eqref{compo}. The obvious
injection of $\fra g$ into $\mathrm L$ induces a direct sum
decomposition
\[
\mathrm L \cong \widetilde {\mathrm L} \oplus \fra g
\]
of chain complexes. Moreover, the Lie algebra perturbation
$\partial_{\Delta}$ on $\mathrm L$ vanishes on the direct summand
$\fra g$ and the other direct summand $\widetilde {\mathrm L}$ is
closed under the operator $\partial_{\Delta}$. Let $\widetilde
{\mathcal L}=\widetilde {\mathrm L}_{\partial_{\Delta}}$; that is to say,
the graded Lie algebra which underlies $\widetilde {\mathcal L}$
coincides with that underlying the kernel $\widetilde {\mathrm L}$
whereas the differential is perturbed via the diagonal map
$\Delta$ of $\Sigm^{\mathrm c}[s \fra g]$. Thus the canonical
projection from $\mathrm L$ to $\fra g$ is also compatible with
the perturbed differential relative to the diagonal map of
$\Sigm^{\mathrm c}[s \fra g]$, and $\widetilde {\mathcal L}$ is
the kernel of the resulting projection $\pi_{\mathcal L}$ from
$\mathcal L$ to $\fra g$. Furthermore, as a chain complex,
$\mathcal L=\mathrm L_{\partial_{\Delta}}$ decomposes as the
direct sum
\[
\mathcal L = \widetilde {\mathcal L} \oplus \nabla_{\mathcal
L}(\fra g),
\]
and $\widetilde {\mathcal L}$ is a differential graded Lie ideal
of $\mathcal L$. Thus the obvious injection $\nabla_{\mathcal
L}\colon \fra g \to \mathcal L$ of $\fra g$ into $\mathcal L$ is a
chain map and the obvious projection $\pi_{\mathcal
L}\colon\mathcal L \to\fra g$ of $\mathcal L$ onto $\fra g$ is a
morphism of differential graded Lie algebras, $\fra g$ being viewed abelian.

For $j \geq 0$, we denote by $\mathcal S^j$ the $j$-th homogeneous
constituent of the symmetric algebra functor $\mathcal S$.

\begin{Lemma}\label{LL}
The homotopy $h_{\Sigm}$ in the contraction \eqref{cont51} induces
a homotopy $h_{\mathcal L}$ such that the data
\begin{equation}
   \left(\fra g
     \begin{CD}
      \null @>{\nabla_{\mathcal L}}>> \null\\[-3.2ex]
      \null @<<{\pi_{\mathcal L}}< \null
     \end{CD}
    \mathcal L\Sigm^{\mathrm c}[s \fra g], h_{\mathcal L}\right)
\label{cont52}
  \end{equation}
constitute a contraction of chain complexes.
\end{Lemma}

\begin{proof}
Consider the perturbed objects
\[
\mathcal L \Sigm^{\mathrm c}[s\fra g]=\mathrm L_{\Delta}[s^{-1}
J\Sigm^{\mathrm c}[s\fra g]],\quad \Omega \Sigm^{\mathrm c}[s\fra
g]= \mathrm T_{\Delta}[s^{-1} J\Sigm^{\mathrm c}[s\fra g]],
\]
the perturbations---beware, not to be confused with the
perturbation $\partial$ defining the sh-Lie algebra structure on
$\fra g$---being induced by the diagonal map of $\Sigm^{\mathrm
c}[s\fra g]$. Relative to the corresponding perturbed
differentials, the projection to the associated graded object
induces an isomorphism
\begin{equation}
\Omega \Sigm^{\mathrm c}[s\fra g] \longrightarrow R \oplus
\mathcal L \Sigm^{\mathrm c}[s\fra g] \oplus \Sigm^2 \mathcal L
\Sigm^{\mathrm c}[s\fra g] \oplus \ldots \oplus \Sigm^{\ell}
\mathcal L \Sigm^{\mathrm c}[s\fra g] \oplus \ldots \label{dir}
\end{equation}
of chain complexes. Furthermore, relative to the direct sum
decomposition \eqref{dir}, for $\ell \geq 1$, the component
\[
\Sigm^{\ell} \mathcal L \Sigm^{\mathrm c}[s\fra g] \longrightarrow
\Sigm^{\ell} \mathcal L \Sigm^{\mathrm c}[s\fra g]
\]
of the homotopy $h_{\Sigm}$ in \eqref{cont51} above is itself a
homotopy and, for $\ell'\ne \ell$, a component of the
kind
\[
\Sigm^{\ell} \mathcal L \Sigm^{\mathrm c}[s\fra g] \longrightarrow
\Sigm^{\ell'} \mathcal L \Sigm^{\mathrm c}[s\fra g],
\]
if non-zero,  is necessarily a cycle (in the corresponding
Hom-complex), since the right-hand side of
\eqref{dir} is a direct sum decomposition of
chain complexes. The component
\[
\mathcal L \Sigm^{\mathrm c}[s\fra g]=\Sigm^{1} \mathcal L
\Sigm^{\mathrm c}[s\fra g] \longrightarrow \Sigm^{1} \mathcal L
\Sigm^{\mathrm c}[s\fra g] =\mathcal L \Sigm^{\mathrm c}[s\fra g]
\]
yields the homotopy $h_{\mathcal L}$ we are looking for.
\end{proof}

We now prove Theorem \ref{lem13} (the sh-Lie algebra perturbation
lemma): Given the contraction \eqref{cont1}, suppose that $\fra g$
comes with a {\em general\/} sh-Lie algebra structure, that is,
let $\partial$ be a {\em general coalgebra perturbation\/} of the
differential $d$ on $\Sigm^{\mathrm c}[s \fra g]$ induced by the
differential on $\fra g$.

The coaugmentation filtration
$\{\mathrm F_n(\Sigm^{\mathrm c}[s \fra g])\}_{(n\geq 0)}$ of
$\Sigm^{\mathrm c}[s \fra g]$ turns $\mathcal L\Sigm^{\mathrm c}[s
\fra g]$ into a filtered differential graded Lie algebra
$\{\mathrm F_n(\mathcal L\Sigm^{\mathrm c}[s \fra g])\}_{(n\geq
0)}$ via
\[
\mathrm F_0(\mathcal L\Sigm^{\mathrm c}[s \fra g]) = 0,\
\mathrm F_n(\mathcal L\Sigm^{\mathrm c}[s \fra g]) =
\mathcal L \mathrm F_n(\Sigm^{\mathrm c}[s \fra g])\quad (n\geq
0),
\]
and we make $\fra g$ into a trivially filtered chain complex
$\{\mathrm F_n(\fra g)\}_{(n\geq 0)}$ via $\mathrm F_0(\fra g)=0$
and $\mathrm F_n(\fra g) = \fra g$ for $n\geq 1$. This turns
\eqref{cont52} into a filtered contraction of chain complexes.
Furthermore, the sh-Lie algebra structure $\partial$ on $\fra g$
(coalgebra perturbation on $\Sigm^{\mathrm c}[s \fra g]$) perturbs
the differential on $\Sigm^{\mathrm c}[s \fra g]$ and hence that
on $\mathcal L\Sigm^{\mathrm c}[s \fra g]$ and, indeed, yields a
Lie algebra perturbation on $\mathcal L\Sigm^{\mathrm c}[s \fra
g]$; we write this perturbation as
\[
\partial_{\mathcal L} \colon \mathcal L\Sigm^{\mathrm c}[s \fra g]
\longrightarrow \mathcal L\Sigm^{\mathrm c}[s \fra g].
\]
Thus perturbing
the loop Lie algebra $\mathcal L\Sigm^{\mathrm
c}[s \fra g]$ on $\Sigm^{\mathrm c}[s \fra g]$
via
$\partial_{\mathcal L}$
carries
the loop Lie algebra 
$\mathcal L\Sigm^{\mathrm c}[s \fra g]$
to the loop Lie algebra
$\mathcal L\Sigm^{\mathrm
c}_{\partial}[s \fra g]$ on $\Sigm^{\mathrm c}_{\partial}[s \fra
g]$.
Application of the ordinary perturbation
lemma (reproduced in \cite{pertlie} as Lemma 5.1) to the Lie
algebra perturbation $\partial_{\mathcal L}$ on $\mathcal
L\Sigm^{\mathrm c}[s \fra g]$ and the filtered contraction of
chain complexes \eqref{cont52} yields the contraction
\begin{equation}
   \left(\fra g
     \begin{CD}
      \null @>{\nabla_{\partial}}>> \null\\[-3.2ex]
      \null @<<{\pi_{\partial}}< \null
     \end{CD}
    \mathcal L\Sigm^{\mathrm c}_{\partial}[s \fra g], h_{\partial}\right)
\label{cont53}
  \end{equation}
of chain complexes. In the special case where the perturbation
$\partial$ arises from an ordinary differential graded Lie algebra
structure on $\fra g$, the morphism $\pi_{\partial}$ is the adjoint of the
resulting Lie algebra twisting cochain $\mathcal C[\fra g] \to
\fra g$ relative to the Lie algebra structure on $\fra g$ and is
therefore a morphism of differential graded Lie algebras relative
to the Lie algebra structure on $\fra g$. Whether or not the
perturbation $\partial$ arises from an ordinary differential
graded Lie algebra structure on $\fra g$, we now combine the
contraction \eqref{cont53} with the original contraction
\eqref{cont1} to the contraction
\begin{equation}
   \left(M
     \begin{CD}
      \null @>{\nabla}>> \null\\[-3.2ex]
      \null @<<{\pi}< \null
     \end{CD}
    \mathcal L\Sigm^{\mathrm c}_{\partial}[s \fra g], h\right)
\label{cont54}
  \end{equation}
of chain complexes where the notation $\nabla$, $\pi$, $h$ is
abused somewhat. More precisely, when the two contractions
\eqref{cont53} and \eqref{cont1} are written as
\begin{equation*}
   \left(M
     \begin{CD}
      \null @>{\nabla_1}>> \null\\[-3.2ex]
      \null @<<{\pi_1}< \null
     \end{CD}
    \fra g, h_1\right),\quad
  \left(\fra g
     \begin{CD}
      \null @>{\nabla_2}>> \null\\[-3.2ex]
      \null @<<{\pi_2}< \null
     \end{CD}
    \mathcal L\Sigm^{\mathrm c}_{\partial}[s \fra g], h_2\right),
  \end{equation*}
the three morphisms in the contraction \eqref{cont54} are given by
\[
\pi = \pi_1 \pi_2,\quad \nabla = \nabla_2 \nabla_1,\quad h =
h_2+\nabla_2 h_1 \pi_2 .
\]
Applying the ordinary Lie algebra perturbation lemma (Lemma
\ref{lem11} above) to the contraction \eqref{cont54} relative to
the differential graded Lie algebra structure on $\mathcal
L=\mathcal L\Sigm^{\mathrm c}_{\partial}[s \fra g]$, we obtain the
perturbation $\pppartial$ on $\Sigm^{\mathrm c}[sM]$, the Lie
algebra twisting cochain
\[
\tau \colon \Sigm_{\ppartial}^{\mathrm c}[sM] \longrightarrow
\mathcal L,
\]
and the asserted contraction \eqref{cont555} of chain complexes,
where we use the notation $\Pi_{\partial}$ and  $H_{\partial}$
rather than the notation $\Pi$ and $H$, respectively, in the
contraction \eqref{cont5} spelled out in the ordinary Lie algebra
perturbation lemma. This completes the proof of Theorem
\ref{lem13}.

\section{Inverting the retraction as an sh-Lie map}
\label{invert}

We return to the situation of the ordinary Lie algebra
perturbation lemma (Lemma \ref{lem11} above). Thus $\fra g$ is now an ordinary
differential graded Lie algebra. Let $\tau$ be the
Lie algebra twisting cochain \eqref{tc1}. The retraction
\[
\Pi \colon \mathcal C[\fra g] \longrightarrow
\Sigm_{\pppartial}^{\mathrm c}[sM]
\]
for the contraction \eqref{cont5} constructed in the last section
of \cite{pertlie} is not in general compatible with the graded
coalgebra structures. As already pointed out, the reason is that 
the notion of homotopy of morphisms of differential graded
cocommutative coalgebras is a subtle concept. We will now explain
how, in the special case where $M$ and $\fra g$ are connected,
the retraction $\Pi$ can be extended to a morphism of sh-Lie algebras, 
that is, to a morphism preserving the appropriate structure.

For intelligibility, we recall the constructions of the retraction
$\Pi$ and  contracting homotopy $H$ in \eqref{cont5} carried out
in \cite{pertlie}: Application of the ordinary perturbation lemma
(reproduced in \cite{pertlie} as Lemma 5.1) to the perturbation
$\partial$ on $\Sigm^{\mathrm c}[s\fra g]$ determined by the
graded Lie algebra structure on $\fra g$ and the induced {\em
filtered\/} contraction
\begin{equation}
   \left(\Sigm^{\mathrm c}[sM]
     \begin{CD}
      \null @>{\Sigm^{\mathrm c}[s\nabla]}>> \null\\[-3.2ex]
      \null @<<{\Sigm^{\mathrm c}[s\pi]}< \null
     \end{CD}
    \Sigm^{\mathrm c}[s\fra g],\Sigm^{\mathrm c}[sh] \right)
\label{2.3}
  \end{equation}
of {\em coaugmented differential graded coalgebras\/}, the
filtrations being the ordinary coaugmentation filtrations, yields
the perturbation $\delta$ of the differential $d^0$ on
$\Sigm^{\mathrm c}[sM]$ and, furthermore, the contraction
\begin{equation}
   \left(\Sigm^{\mathrm c}_{\delta}[sM]
     \begin{CD}
      \null @>{\nnabla}>> \null\\[-3.2ex]
      \null @<<{\widetilde \Pi}< \null
     \end{CD}
    \mathcal C[\fra g],\HH \right)
\label{2.6a}
  \end{equation}
of chain complexes. Moreover, the composite
\begin{equation}
\begin{CD}
\Phi\colon\Sigm^{\mathrm c}_{\pppartial}[sM]
@>{\overline \tau}>> \mathcal C[\fra g] @>{\widetilde \Pi}>>
\Sigm^{\mathrm c}_{\delta}[sM]
\end{CD}
\label{comp2}
\end{equation}
is an isomorphism of chain complexes, and the morphisms
\begin{align}
\Pi &= \Phi^{-1}\widetilde \Pi \colon \mathcal C[\fra g]
\longrightarrow \Sigm^{\mathrm c}_{\pppartial}[sM],
\label{Pipert1}
\\
H &= \HH- \HH \overline \tau\,\Pi \colon \mathcal C[\fra g]
\longrightarrow \mathcal C[\fra g] \label{Hpert1}
\end{align}
complete the construction of the contraction \eqref{cont5}.

In general, none of the morphisms $\delta$, $\nnabla$, $\widetilde
\Pi$, $\Pi$, $\HH$, $H$ is compatible with the coalgebra
structures. The isomorphism of chain complexes $\Phi$ admits an
explicit description in terms of the data as a {\em perturbation
of the identity\/}  and so does its inverse; details have been
given in the last section of \cite{pertlie}.

Consider the universal loop Lie algebra twisting cochain
\begin{equation}
t_{\mathcal L} \colon \Sigm_{\pppartial}^{\mathrm c}[sM]
\longrightarrow \mathcal L\Sigm_{\pppartial}^{\mathrm c}[sM] .
\label{tc5}
\end{equation}
We recall that $M$ to be connected means that
$M$ is concentrated either in positive or in negative degrees;
in particular, the degree zero constituent of $M$ is zero.

\begin{Lemma}
\label{lem72} Suppose that $M$ is connected.
The recursive construction
\begin{equation}
\vartheta = t_{\mathcal L} \Pi +\tfrac[\vartheta,\vartheta] H \colon
\mathcal C[\fra g] \longrightarrow \mathcal
L\Sigm_{\pppartial}^{\mathrm c}[sM] \label{lieloop1}
\end{equation}
yields a Lie algebra twisting cochain $\vartheta\colon \mathcal
C[\fra g] \longrightarrow \mathcal L\Sigm_{\pppartial}^{\mathrm
c}[sM]$ such that
\begin{equation}
\vartheta \overline \tau =t_{\mathcal L} \colon
\Sigm_{\pppartial}^{\mathrm c}[sM] \longrightarrow \mathcal
L\Sigm_{\pppartial}^{\mathrm c}[sM] . \label{tc3}
\end{equation}
\end{Lemma}

\begin{proof} The construction \eqref{lieloop1} being recursive means that
\[
\vartheta=  \vartheta_1 + \vartheta_2 + \ldots
\]
where $\vartheta_1= t_{\mathcal L} \Pi$, $\vartheta_2=
\tfrac 12[\vartheta_1,\vartheta_1]H $,
$\vartheta_3=[\vartheta_1,\vartheta_2]H$, etc. 
The connectedness hypothesis entails the convergence, which is naive.
We leave the
details as an exercise. \end{proof}

\medskip
\noindent {\bf Complement I to Lemma \ref{lem11}}. {\sl In view of
the identity  \/} \eqref{tc3}, {\sl it is manifest that the
composite
\begin{equation*}
\begin{CD}
\Sigm_{\pppartial}^{\mathrm c}[sM] @>{\overline \tau}>>
\mathcal C[\fra g]
@>{\vartheta}>>\mathcal L\Sigm_{\pppartial}^{\mathrm c}[sM]
\end{CD}
\end{equation*}
coincides with the universal loop Lie algebra twisting cochain
\eqref{tc5}. In this sense, $\vartheta$ yields  an sh-retraction
for the sh-morphism from $(M,\pppartial)$ to $\fra g$ given by
$\overline \tau$.}

\smallskip

To explain in which sense the other composite
\begin{equation}
\begin{CD}
\fra g @>{\vartheta}>> (M,\pppartial) @>{\overline\tau}>>
\fra g
\end{CD}
\label{comp}
\end{equation}
of these morphisms is homotopic to the identity, we need some more
preparation.

Let $C$ be a coaugmented differential graded coalgebra and $A$ an
augmented differential graded algebra. Recall that, given two
ordinary twisting cochains $\tau_1,\tau_2\colon C \to A$, a {\em
homotopy\/}
$
h\colon \tau_1 \simeq \tau_2
$
of twisting cochains is a homogeneous morphism
\begin{equation}
h \colon C \longrightarrow A \label{tc7}
\end{equation}
of degree zero such that $\varepsilon h \eta = \varepsilon \eta$
and
\begin{equation}
Dh =\tau_1 \cup h - h \cup \tau_2 \in \mathrm{Hom}(C,A).
\label{tc10}
\end{equation}
Such a homotopy $h\colon \tau_1 \simeq \tau_2$ of twisting
cochains is well known to induce a chain homotopy
\begin{equation}
\overline h \colon C \longrightarrow BA \label{tc11}
\end{equation}
between the adjoints $\overline \tau_1,\overline \tau_2 \colon C
\longrightarrow BA$ into the reduced bar construction $BA$ on $A$, and the
homotopy $\overline h$ is compatible with the coalgebra
structures.

Recall that the augmented
differential graded algebra $A$ is {\em complete\/} when
the canonical morphism of differential graded algebras from $A$ to
$\lim A\big/(IA)^n$ is an isomorphism; here $IA$ refers to the 
augmentation ideal as usual.

\begin{Lemma}
\label{lem6} Suppose the following data are given:

\noindent
--- coaugmented differential graded coalgebras $B$ and
$C$;

\noindent
--- a contraction
\Nsddata C{\nabla}{\pi}Bh of chain complexes, $\nabla$ being a
morphism of coaugmented differential graded coalgebras;

\noindent
--- an augmented differential graded algebra $A$;

\noindent
--- twisting cochains $t_1,t_2 \colon C \to A$;

\noindent
--- a homotopy $h^B\colon B \to A$ of twisting cochains $h^B\colon
t_1\nabla \simeq t_2\nabla$, so that
\begin{equation}
D(h^B)=(t_1 \nabla)\cup h^B - h^B\cup(t_2 \nabla). \label{tc8}
\end{equation}
Suppose that the augmented differential graded
algebra  $A$ is complete. Then the recursive formula
\begin{equation}
h^C=h^B\pi - (t_1 \cup h^C - h^C\cup t_2)h \label{tc9}
\end{equation}
yields a homotopy $h^C\colon C \to A$ of twisting cochains
$h^C\colon t_1\simeq t_2$ such that $h^C\nabla = h^B$.
\end{Lemma}

The formula \eqref{tc9} being recursive means that
\[
h^C=  \varepsilon \eta + h_1 + h_2 + \ldots
\]
where $h_1=h^B\pi - (t_1 -t_2)h$, $h_2= - (t_1 \cup h_1 - h_1\cup
t_2)h$, etc.

\begin{proof} The identity $h^C\nabla = h^B$ is obvious and,
since $t_1$ and $t_2$ are ordinary twisting cochains, the morphism
$t_1 \cup h^C - h^C\cup t_2$ is (easily seen to be) a cycle.
Furthermore, since $\nabla$ is compatible with the coalgebra
structures,
\begin{align*}
(t_1 \cup h^C - h^C\cup t_2)\nabla \pi &=((t_1\nabla) \cup
(h^C\nabla) - (h^C\nabla)\cup (t_2\nabla)) \pi \\&= ((t_1
\nabla)\cup h^B - h^B\cup(t_2 \nabla))\pi .
\end{align*}
Consequently
\begin{align*}
Dh^C &= (D(h^B))\pi + (t_1 \cup h^C - h^C\cup t_2)Dh
\\&=(t_1 \nabla)\cup h^B - h^B\cup(t_2 \nabla)\pi +
(t_1 \cup h^C - h^C\cup t_2)- (t_1 \cup h^C - h^C\cup t_2)\nabla
\pi
\\
&= t_1 \cup h^C - h^C\cup t_2
\end{align*}
as asserted. \end{proof}

Let $(\fra h_1,
\partial_1)$ and $(\fra h_2,\partial_2)$ be two sh-Lie algebras
and let
\[
\vartheta_1,\vartheta_2\colon \Sigm^{\mathrm c}_{\partial_1}[s\fra
h_1] \longrightarrow \mathcal L \Sigm^{\mathrm
c}_{\partial_2}[s\fra h_2]
\]
be two Lie algebra twisting cochains, that is, generalized
sh-morphisms or generalized sh-Lie maps from $(\fra h_1,
\partial_1)$ to $(\fra h_2, \partial_2)$.
We define a {\em homotopy of generalized sh-morphisms\/} or {\em
homotopy of generalized sh-Lie maps\/} from $\vartheta_1$ to
$\vartheta_2$ to be a homotopy
\begin{equation}
h\colon \Sigm^{\mathrm c}_{\partial_1}[s\fra h_1] \longrightarrow
\mathrm U\mathcal L \Sigm^{\mathrm c}_{\partial_2}[s\fra h_2]
=\Omega\Sigm^{\mathrm c}_{\partial_2}[s\fra h_2] \label{htp}
\end{equation}
of ordinary twisting cochains $h\colon \vartheta_1 \simeq
\vartheta_2$. Here and below we identify a Lie algebra twisting
cochain with the corresponding ordinary twisting cochain having
values in the corresponding universal algebra.

\begin{Remark} Write $\mathcal L= \mathcal L \Sigm^{\mathrm
c}_{\partial_2}[s\fra h_2]$. In view of the definitions, the
adjoint of a homotopy $h$ of the kind {\rm \eqref{htp}} takes the
form $ \overline h\colon \Sigm^{\mathrm c}_{\partial_1}[s\fra h_1]
\longrightarrow \mathrm B \mathrm U\mathcal L =\mathrm B \Omega
\Sigm^{\mathrm c}_{\partial_2}[s\fra h_2]$, whence the values of
the adjoint $\overline h$ of the homotopy \eqref{htp} necessarily
lie in the coaugmented differential graded coalgebra $\mathrm B
\mathrm U\mathcal L$ rather than in the coaugmented differential
graded cocommutative coalgebra $\mathcal C[\mathcal L]$, viewed as
a subcoalgebra of $\mathrm B \mathrm U\mathcal L$ via the
canonical injection
\begin{equation}
\mathcal C[\mathcal L] \longrightarrow \mathrm B \mathrm U\mathcal
L.\label{inject}
\end{equation}
The injection \eqref{inject}, in turn, is well known to be a
quasi-isomorphism, though.

Historically, the injection \eqref{inject} has played a major role
for the development of Lie algebra cohomology, cf. e.~g. {\rm
\cite{cartanei} (Ch. XIII, Theorem 7.1)} for the special case of
an ordinary (ungraded) Lie algebra. From the point of view of
sh-Lie algebras, $\mathcal C[\mathcal L]$ would be the correct
target for the adjoint  of a homotopy
of the kind
\eqref{htp}. To arrive at an adjoint
having values in $\mathcal C[\mathcal L]$, one would
have to require that the values of a homotopy of twisting
cochains of the kind \eqref{htp} lie in $\mathcal L$ rather than in $\mathrm
U[\mathcal L]= \Omega\Sigm^{\mathrm c}_{\partial_2}[s\fra h_2]$.
Such a requirement would lead to inconsistencies, though: The
requirement that a homotopy of the kind $\overline h$ be
compatible with coalgebra structures forces a condition of the
kind \eqref{tc10}; this condition, in turn, necessarily involves
the multiplication map in the universal algebra $\mathrm U\mathcal
L=\Omega \Sigm^{\mathrm c}_{\partial_2}[s\fra h_2]$ of the
corresponding differential graded Lie algebra $\mathcal L$ (rather
than just the graded Lie algebra structure of $\mathcal L$)
and hence cannot be phrased
merely in terms of the graded Lie algebra structure
alone, whence the values of the
homotopy \eqref{htp} cannot in general lie in $\mathcal
L$. Thus, strictly speaking, the notion of {\em homotopy leaves
the world of sh-Lie algebras\/}. Again this observation reflects
the fact that the notion of homotopy of morphisms of differential
graded cocommutative coalgebras is a subtle concept.

Nevertheless, a cure is provided by an appropriate higher
homotopies construction:  A differential graded coalgebra of the
kind $\mathrm B \mathrm U\mathcal L=\mathrm B \Omega
\Sigm^{\mathrm c}_{\partial_2}[s\fra h_2]$ is a {\em
quasi-commuted\/} coalgebra, cf. {\rm \cite{husmosta} (p.~175)};
moreover, in the category {\rm DCSH, cf. \cite{gugenmun}}, the
injection {\rm \eqref{inject}} is an isomorphism (preserving the
diagonal maps), and the diagonal map of $\mathrm B \mathrm
U\mathcal L=\mathrm B \Omega \Sigm^{\mathrm c}_{\partial_2}[s\fra
h_2]$ is a morphism in the category. Thus, suitably rephrased, the
notion of homotopy will stay within the world of sh-Lie algebras.
The exploration of categories of the kind {\rm DCSH} was 
prompted by {\rm \cite{halpstas}}.
\end{Remark}

We will now exploit Lemma \ref{lem6} in the following manner: 
Suppose that $M$ and $\fra g$ are connected.
Let
$B=\Sigm_{\ppartial}^{\mathrm c}[sM]$, $C=\mathcal C[\fra g]$,
consider the contraction \eqref{cont5}, 
let $A=\mathrm U\mathcal L\mathcal C[\fra g]=\Omega \mathcal
C[\fra g]$---notice that $A$ is connected in the sense
that $A_0$ is a copy of the ground ring---, and let
\begin{align*}
t_1&=\mathcal L(\overline\tau)\vartheta \colon \mathcal C[\fra g]
\longrightarrow \mathcal L\mathcal C[\fra g],
\\
t_2&=t_{\mathcal L} \colon \mathcal C[\fra g] \longrightarrow
\mathcal L\mathcal C[\fra g],
\\
h^B&=\varepsilon \eta.
\end{align*}
By construction,
\[
t_1 \overline \tau = t_2 \overline\tau \colon
\Sigm_{\ppartial}^{\mathrm c}[sM] \longrightarrow \mathcal
L\mathcal C[\fra g],
\]
and Lemma \ref{lem6} applies. These observations establish the
following.

\medskip
\noindent {\bf Complement II to Lemma \ref{lem11}}. {\sl 
Suppose that $\fra g$ is connected.
The
homotopy 
\[
 h^C\colon \mathcal C[\fra g] \longrightarrow \mathrm
U\mathcal L\mathcal C[\fra g] =\Omega\mathcal C[\fra g]
\]
of
twisting cochains $h^C\colon t_1\simeq t_2$ given by {\rm
\eqref{tc9}} yields a homotopy between the composite\/} {\rm
\eqref{comp}} {\sl and the identity of $\fra g$, all objects and
morphisms in sight being viewed as sh-objects and sh-morphisms.}

\medskip

Constructions of the same kind yield an explicit sh-inverse for
\eqref{tc2} as a twisting cochain of the kind
\[
\mathcal C[\mathcal L \Sigm_{\partial}^{\mathrm c}[s\fra g] ]
\longrightarrow \mathcal L \Sigm_{\ppartial}^{\mathrm c}[sM]
\]
as well, $M$ and $\fra g$
still being supposed to be connected. 
We spare the reader and ourselves these added troubles
here.

\section{The proof of the theorem in the introduction}
\label{proofthm}

Let $\partial$ be an sh-Lie algebra structure on $\fra g$, and let
$\ppartial$ be the coalgebra perturbation on $\Sigm^{\mathrm c}[sM]$
and
\[
\tau\colon\Sigm_{\ppartial}^{\mathrm
c}[sM]\longrightarrow \mathcal L\Sigm^{\mathrm
c}_{\partial}[s \fra g]
\]
the
Lie algebra twisting cochain \eqref{tc2} given in
the sh-Lie algebra perturbation lemma (Theorem \ref{lem13}
above). The theorem in the introduction
comes down to the observation that, with the notation of the
previous section, both the adjoint
\begin{equation}
\overline \tau\colon\Sigm_{\ppartial}^{\mathrm
c}[sM]\longrightarrow \mathcal C[\mathcal L\Sigm^{\mathrm
c}_{\partial}[s \fra g]] \label{adjointt1}
\end{equation}
of $\tau$  and the adjoint
\begin{equation}
\overline {t_{\mathcal L}}\colon \Sigm^{\mathrm c}_{\partial}[s
\fra g] \longrightarrow \mathcal C[\mathcal L\Sigm^{\mathrm
c}_{\partial}[s \fra g]] \label{adjointt2}
\end{equation}
of the universal loop Lie algebra twisting  cochain $ t_{\mathcal
L}\colon \Sigm^{\mathrm c}_{\partial}[s \fra g] \longrightarrow
\mathcal L\Sigm^{\mathrm c}_{\partial}[s \fra g] $
yield
sh-equiva\-len\-ces. Under appropriate connectivity hypotheses,
constructions similar to those spelled out in
the previous section yield explicit sh-inverses for
\eqref{adjointt1} and \eqref{adjointt2}.

\end{document}